\newfont{\bcb}{msbm10}
\newfont{\matb}{cmbx10}
\newfont{\got}{eufm10}

\documentclass[12pt]{article}
\usepackage[cp1250]{inputenc}

\usepackage{amsmath,amsthm}
\usepackage{amssymb,latexsym}
\usepackage{enumerate}

\newcommand{\matR}{\mathbb{R}}

\newcommand{\matN}{\mathbb{N}}

\begin{document}

\title{On division of quasianalytic function germs}

\author{by Krzysztof Jan Nowak}




\date{}

\maketitle


\begin{abstract}
In this paper, we establish the following criterion for
divisibility in the local ring $\mathcal{Q}_{n}$ of those
quasianalytic function germs at $0 \in \matR^{n}$ which are
definable in a polynomially bounded structure. A sufficient (and
necessary) condition for the divisibility of two function germs in
$\mathcal{Q}_{n}$ is that of their Taylor series at $0 \in
\matR^{n}$ in the formal power series ring.
\end{abstract}


Fix a polynomially bounded, o-minimal structure $\mathcal{R}$ on
the field of reals and denote by $\mathcal{Q}_{n}$ the ring of
those quasianalytic function germs at $0 \in \matR^{n}$ which are
definable in $\mathcal{R}$. The main purpose of this paper is to
establish the following

\vspace{2ex}

{\bf Criterion for Divisibility.}
\begin{em}
Let $f,g \in \mathcal{Q}_{n}$ be two quasianalytic germs at $0 \in
\matR^{n}$ and $\widehat{f}, \widehat{g} \in \matR [[x]]$,
$x=(x_{1},\ldots,x_{n})$, be their Taylor series at $0 \in
\matR^{n}$. If $\; \widehat{f} \in \widehat{g} \cdot \matR [[x]]$,
then $\; f \in g \cdot \mathcal{Q}_{n}$. In other words,
\end{em}
$$ g \cdot \matR [[x]] \, \cap \, \mathcal{Q}_{n} = g \cdot \mathcal{Q}_{n}. $$


Before proceeding with the proof, we make a comment. Our criterion
is much stronger than Bierstone--Milman's quasianalytic division
theorem (cf.~\cite{BM2}, Theorem~6.4), and is not a local version
of the latter, which says that a smooth function $f$ on an open
subset $U$ of $\mathbb{R}^{n}$ is divisible by a quasianalytic
function $g$ on $U$, if it is formally divisible by $g$ (i.e.\ at
all points $a \in U$, the Taylor series of $f$ at $a$ is divisible
by the Taylor series of $g$ at $a$). (Their proof followed
Atiyah's proof of the division theorem of H\"{o}rmander and
\L{}ojasiewicz.) Whereas we assume only the divisibility of the
Taylor series at one point (zero), they assume the divisibility of
the Taylor series at all points of a given open subset $U$. Our
criterion implies, therefore, that the divisibility of the Taylor
series propagates from one point into its vicinity. Such a result
is very subtle in quasianalytic geometry because quasianalytic
local rings enjoy neither good algebraic properties nor coherence.

It should be emphasized here that the problem whether
quasianalytic local rings (larger than the analytic ones) are
Noetherian remains open as yet (cf.~\cite{Elkh}). It is only known
that quasianalytic local rings coming from Denjoy--Carleman
classes do not admit the Weierstrass division theorem
(cf.~\cite{Chil}) and even the Weierstrass preparation theorem
(cf.~\cite{ABBNZ}). Neither do more general quasianalytic local
rings defined by imposing certain natural axioms (stability by
composition, implicit function theorem and monomial division); see
\cite{Elkh-S} for Weierstrass division and \cite{Par-Rol} for
Weierstrass preparation.

Finally, let us mention that the quasianalytic division theorem
and formal divisibility was also studied in
papers~\cite{Bew,Now25}.

\vspace{2ex}

{\em Proof.} Since the division problem is local, we may assume
that the structure $\mathcal{R}$ is the real field with restricted
Q-analytic functions determined by the local rings
$\mathcal{Q}_{n}$, $n \in \matN$. Notice that Q-analytic functions
fulfill the above-mentioned axioms, which ensure desingularization
algorithms, upon which the geometry of quasianalytic structures
relies (cf.~\cite{Rol-Spei-Wil,Now1,Now2,Now5}). We thus have at
our disposal a powerful tool, transformation to a normal crossing
by blowing up (cf.~\cite{BM2,Rol-Spei-Wil}), and even canonical
one, along global smooth centers (cf.~\cite{BM1}), which will be
applied afterwards.


The structure $\mathcal{R}$ admits, of course, smooth cell
decompositions and quasianalytic stratifications
(cf.~\cite{Rol-Spei-Wil,Now1,Now2}). Our proof of the criterion
relies on a stabilization theorem from our paper~\cite{Now4} to
the effect that the sets $\mathcal{C}_{k}(f)$ of points near which
a given definable function $f$ is of class $\mathcal{C}^{k}$, $k
\in \matN \cup \{ \infty \}$, stabilizes. This theorem was
achieved, in turn, by means of some stabilization effects linked
to Gateaux differentiability and to formally composite functions.
The latter relied also on a quasianalytic version of the composite
function theorem from our paper~\cite{Now3}, where we adapted the
approach developed by Bierstone--Milman--Pawłucki~\cite{BMP}.


For simplicity, we shall use the same letters to denote functions
and the induced germs. By this abuse of notation, we regard $f,g$
as functions Q-analytic in a common neighbourhood $U$ of $0 \in
\matR^{n}$. Let $\mathrm{T}_{a}f(x)$ stand for the Taylor series
of a function $f$ at a point $a$.

\vspace{1ex}

We first establish the following

\vspace{2ex}

{\bf Claim.}
\begin{em}
With the above notation, if $\; \widehat{f} \in \widehat{g} \cdot
\matR [[x]]$, then the quotient $f/g$ extends to a unique
continuous function $\varphi$ definable in a neighbourhood of
zero.
\end{em}

\vspace{2ex}

For the proof, observe that there exists a transformation $\sigma:
\widetilde{U} \longrightarrow U$ to a normal crossing by global
blowings-up such that, in a suitable local coordinate system $y$
near any point $b \in \widetilde{U}$, we have
$$ f^{\sigma}(y) := f \circ \sigma (y) = \mbox{unit} \cdot y^{\beta}, \ \
   g^{\sigma}(y) := g \circ \sigma (y) = \mbox{unit} \cdot y^{\gamma}, $$
where $\beta,\gamma \in \matN^{n}$ (being dependent on the point
$b$), and $\beta \geq \gamma$ or $\beta \leq \gamma$; here $\leq$
stands for the partial ordering determined by divisibility
relation. By assumption,
$$ H(x) := f(x)/g(x) = \mathrm{T}_{0}f(x)/\mathrm{T}_{0}g(x) \in
   \matR[[ x]]; $$
here we identify quasianalytic germs at zero with their Taylor
series. Clearly,
$$ \mathrm{T}_{b} \, (f^{\sigma}(y)/g^{\sigma}(y)) =
   \mathrm{T}_{b} \, f^{\sigma}(y)/ \mathrm{T}_{b} \, g^{\sigma}(y) =
   H \circ \mathrm{T}_{b} \, \sigma(y), $$
whence $\beta \geq \gamma$ for any point $b \in \sigma^{-1}(0)$ on
the fiber. Consequently, the quotient $f^{\sigma}/g^{\sigma}$ is
near any point $b \in \sigma^{-1}(0)$ of the form: $\mbox{unit}
\cdot y^{\beta - \gamma}$, and thus extends to a unique Q-analytic
function $\psi$ in the vicinity of the fibre $\sigma^{-1}(0)$ on
$\widetilde{U}$.


We show, by {\em reductio ad absurdum}, that the function $\psi$
is constant on the fibres of the transformation $\sigma$ in a
neighbourhood of $0 \in \matR^{n}$. Suppose the contrary and put
$$ E := \{ (x,y,z) \in U \times \widetilde{U} \times
   \widetilde{U}: \ \sigma(y) = \sigma(z)=x, \ \psi(y) \neq \psi(z)
   \}. $$
Then, since the transformation $\sigma$ is a proper mapping, the
closure of the set $E$ contains a point $(0,b_{0},c_{0}) \in U
\times \widetilde{U} \times \widetilde{U}$ with $\sigma(b_{0})
=\sigma(c_{0})= 0$. By the arc selection lemma along with the
quasianalytic version of Puiseux's theorem (cf.~\cite{Now4},
Section~2), there exist three Q-analytic arcs
$$ a(t): (-1,1) \longrightarrow U, \ \ b(t),c(t): (-1,1)
   \longrightarrow \widetilde{U} $$
such that $a(0)=0$, $b(0)=b_{0}$, $c(0)=c_{0}$ and
$$ \{ (a(t),b(t),c(t)): \; t \in (0,1) \} \subset E. $$
Then
$$ (\sigma \circ b)(t) = (\sigma \circ c)(t) = a(t) \ \ \mbox{ and } \ \
    (\psi \circ b)(t) \neq (\psi \circ c)(t) $$
for all $t \in (0,1)$. Clearly,
$$ H \circ \mathrm{T}_{b_{0}} \, \sigma = \mathrm{T}_{b_{0}} \, \psi \
   \ \mbox{ and } \ \
   H \circ \mathrm{T}_{c_{0}} \, \sigma = \mathrm{T}_{c_{0}} \, \psi.
   $$
We thus get
$$ \mathrm{T}_{0}\, (\psi \circ b) = \mathrm{T}_{b_{0}} \, \psi
   \circ \mathrm{T}_{0}\, b = H \circ \mathrm{T}_{b_{0}} \, \sigma
   \circ \mathrm{T}_{0}\, b = H \circ \mathrm{T}_{0}\, (\sigma
   \circ b) = H \circ \mathrm{T}_{0}\, a, $$
and, similarly, $\mathrm{T}_{0}\, (\psi \circ b) = H \circ
\mathrm{T}_{0}\, a$. Hence $\mathrm{T}_{0}\, (\psi \circ b) =
\mathrm{T}_{0}\, (\psi \circ c)$, and by quasianalyticity, $\psi
\circ b \equiv \psi \circ c$, which is a contradiction.


Therefore, since $\sigma$ is a proper mapping, $\psi$ descends to
a continuous function $\varphi$ in a neighbourhood of $0 \in
\matR^{n}$. Obviously, $\varphi$ is a continuous extension of the
quotient $f/g$, which proves the claim.

\vspace{2ex}

Similarly, it can be immediately deduced from the above claim that
each partial derivative $D^{\alpha}\, f/g$, $\alpha \in
\matN^{n}$, extends to a unique continuous function definable in a
neighbourhood $U_{\alpha}$ of $0 \in \matR^{n}$, too, because
$D^{\alpha}\, f/g$ is a quotient of two quasianalytic germs such
that the quotient of their Taylor series coincides with
$D^{\alpha}\, H$. Notice that the partial derivatives (up to a
given order $k$) are defined off a nowhere dense singular locus of
the quotient $h$ (the zero locus of the denominator $g$) and
extend continuously onto a neighbourhood $U_{k}$ of $0 \in
\matR^{n}$. By the good direction lemma (cf.~\cite{Dries},
Chap.~7, \S~4), one can assume that the lines parallel to the axes
meet the singular locus in finitely many points. Hence and by de
l'Hospital's rule, the (unique) continuous extensions of the
partial derivatives are everywhere the partial derivatives of
$\varphi$. Therefore, the extension $\varphi$ has on $U_{k}$
continuous partial derivatives up to order $k$, and thus is
$\mathcal{C}^{k}$ on $U_{k}$.


Therefore, for each $k \in \matN$, the extension $\varphi$ of the
quotient $f/g$ is of class $\mathcal{C}^{k}$ in a (possibly
smaller and smaller) neighbourhood of $0 \in \matR^{n}$. A priori,
we are not able to find a common neighbourhood for all the above
partial derivatives. However, it follows from our stabilization
theorem (cf.~\cite{Now4}, Section~2) that the set of points where
the function $\varphi$ is smooth (i.e.\ of class
$\mathcal{C}^{\infty}$) coincides with the set of points where it
is of class $\mathcal{C}^{k}$ for $k$ large enough. Consequently,
$\varphi$ is a smooth function definable in a neighbourhood of $0
\in \matR^{n}$, and thus it is an extension of the quotient $f/g$
we are looking for. This completes the proof of the criterion.

\vspace{2ex}

{\bf Remark.} It is most plausible that the quasianalytic local
rings $\mathcal{Q}_{n}$ are not Noetherian. It would be desirable
to find a non-principal ideal $I$ of $\mathcal{Q}_{n}$ such that
$$ I \varsubsetneq I \cdot \mathbb{R}[[ x ]] \cap \mathcal{Q}_{n}, $$
which would confirm, of course, this conjecture. We do not know if
the analogous criterion for divisibility will hold under certain
conditions imposed on quasianalytic Denjoy--Carleman classes.

\vspace{4ex}

\begin{small}
\begin{sc}
Institute of Mathematics

Faculty of Mathematics and Computer Science

Jagiellonian University

ul.~Profesora \L{}ojasiewicza 6

30-348 Krak\'{o}w, Poland

{\em e-mail address: nowak@im.uj.edu.pl}
\end{sc}
\end{small}

\end{document}